\documentclass[11pt]{article}
\usepackage{subfigure}
\usepackage{url}
\usepackage{amsmath, amssymb, amsthm, color}
\usepackage{graphicx}

\newtheorem{proposition}{Proposition}
\newtheorem{theorem}{Theorem}

\newcommand{\bx}{\mbox{\boldmath $x$}}

\newcommand{\by}{\mbox{\boldmath $y$}}
\newcommand{\bz}{\mbox{\boldmath $z$}}

\newcommand{\bb}{\mbox{\boldmath $b$}}
\newcommand{\bp}{\mbox{\boldmath $p$}}
\newcommand{\bN}{\mbox{\boldmath $N$}}

\newcommand{\be}{\mbox{\boldmath $e$}}
\newcommand{\bxi}{\mbox{\boldmath $\xi$}}
\newcommand{\bmu}{\mbox{\boldmath $\mu$}}
\newcommand{\balpha}{\mbox{\boldmath $\alpha$}}
\newcommand{\bbeta}{\mbox{\boldmath $\beta$}}
\newcommand{\ba}{\mbox{\boldmath $a$}}
\newcommand{\br}{\mbox{\boldmath $r$}}

\long\def\omitthis#1{}

\usepackage{float}
\floatstyle{ruled}
\newfloat{algorithm}{htp}{lop}
\floatname{algorithm}{Algorithm}
\begin{document}
\title{Likelihood Robust Optimization for Data-driven Problems}

\author{Zizhuo Wang \thanks{Department of Industrial and Systems Engineering, University of
Minnesota, Minneapolis, MN, 55455. {\tt Email: zwang@umn.edu}.} \and
Peter W. Glynn\thanks{Department of Management Science and
Engineering, Stanford University, Stanford, CA 94305, {\tt
Email:glynn@stanford.edu.}}\and Yinyu Ye\thanks{Department of
Management Science and Engineering, Stanford University, Stanford,
CA 94305, {\tt Email: yinyu-ye@stanford.edu. }}}

\maketitle

\begin{abstract}
We consider optimal decision-making problems in an uncertain
environment. In particular, we consider the case in which the
distribution of the input is unknown, yet there is abundant
historical data drawn from the distribution. In this paper, we
propose a new type of distributionally robust optimization model
called the {\it likelihood robust optimization} (LRO) model for this
class of problems. In contrast to previous work on distributionally
robust optimization that focuses on certain parameters (e.g., mean,
variance, etc.) of the input distribution, we exploit the historical
data and define the accessible distribution set to contain only
those distributions that make the observed data achieve a certain
level of likelihood. Then we formulate the targeting problem as one
of optimizing the expected value of the objective function under the
worst-case distribution in that set. Our model avoids the
over-conservativeness of some prior robust approaches by ruling out
unrealistic distributions while maintaining robustness of the
solution for any statistically likely outcomes. We present
statistical analyses of our model using Bayesian statistics and
empirical likelihood theory. Specifically, we prove the asymptotic
behavior of our distribution set and establish the relationship
between our model and other distributionally robust models. To test
the performance of our model, we apply it to the newsvendor problem
and the portfolio selection problem. The test results show that the
solutions of our model indeed have desirable performance.
\end{abstract}

\section{Introduction}
\label{sec:intro} The study of decision making problems in uncertain
environments has been a main focus in the operations research
community for decades. In such problems, one has a certain objective
function to optimize, however, the objective function depends not
only on the decision variables, but also on some unknown parameters.
Such situations are ubiquitous in practice. For example, in an
inventory management problem, the inventory cost is influenced by
both the inventory decisions and the random demands. Similarly, in a
portfolio selection problem, the realized return is determined by
both the choice of the portfolio and the random market fluctuations.

One solution method to such problems is stochastic optimization. In
this approach, one assumes the knowledge of the distribution of the
unknown parameters and chooses the decision that optimizes the
expected value of the objective function. If the knowledge of the
distribution is exact, then this approach is a precise
characterization for a risk-neutral decision maker. Much research
has been conducted on this topic; we refer the readers to Shapiro et
al. \cite{shapiro} for a comprehensive review on the topic of
stochastic optimization.

However, there are several drawbacks to the stochastic optimization
approach. First, although stochastic optimization can frequently be
formulated as a convex program, in order to solve it, one often has
to resort to the Monte Carlo method, which can be computationally
challenging. More importantly, due to the limitation of knowledge,
the distribution of the uncertain parameters is rarely known in
practice to a precise level. Even if enough data have been gathered
in the past to perform statistical analyses for the distribution,
the analyses are often based on assumptions (e.g., independence of
the observations, or stationarity of the sequence) that are only
approximations of the reality. In addition, many decision makers in
practice are not risk-neutral. They tend to be risk-averse. A
solution approach that can guard them from adverse scenarios is of
great practical interest.

One such approach was proposed by Scarf \cite{scarf} in a newsvendor
context and has been studied extensively in the past decade. It is
called the {\it distributionally robust optimization} (DRO)
approach. In the DRO approach, one considers a {\it set} of
distributions for the uncertain parameters and optimizes the
worst-case expected value of the objective function among the
distributions in that set. Studies have been done by choosing
different distribution sets. Among them, most choose the
distribution set to contain those distributions with a fixed mean
and variance. For example, Scarf \cite{scarf} shows that a
closed-form solution can be obtained in the newsvendor problem
context when such a distribution set is chosen. Other earlier works
include Gallego and Moon \cite{gallego}, Dupa\^{c}ov\'{a}
\cite{dupacova} and \^{Z}\'{a}\^{c}kov\'{a} \cite{zackova}. The same
form of the distribution set is also used in Calafiore and El Ghaoui
\cite{calafiore}, Yue et al. \cite{yue}, Zhu et al. \cite{zhu} and
Popescu \cite{popescu} in which a linear-chance-constrained problem,
a minimax regret objective and a portfolio optimization problem are
considered, respectively. Other distribution sets beyond the mean
and variance have also been proposed in the literature. For example,
Delage and Ye \cite{delage} propose a more general framework with a
distribution set formed by moment constraints. A review of the
recent developments can be found in Delage \cite{erickthesis} and
Shapiro and Klegwegt \cite{shapiro2}.\footnote{We note that there is
also a vast literature on robust optimization where the worst-case
parameter is chosen for each decision made. However, the robust
optimization is based on a slightly different philosophy than the
distributionally robust optimization and is usually more
conservative. It can also be viewed as a special case of the
distributionally robust optimization where the distribution set only
contains singleton distributions. In view of this, we choose not to
include a detailed discussion of this literature in the main text
and refer the readers to Ben-Tal et al. \cite{bentalbook} and
Bertsimas et al. \cite{brown} for comprehensive reviews.}

Although the mean-variance DRO approach is intuitive and is
tractable under certain conditions, it is unsatisfactory from at
least two aspects. First, when constructing the distribution set in
such an approach, one only uses the moment information in the sample
data, while all the other information is ignored. This procedure may
discard important information in the data set. For example, a set of
data drawn from an exponential distribution with $\lambda = 1/50$
will have similar mean and variance as a set of data drawn from a
normal distribution with $\mu =\sigma = 50$. In the mean-variance
DRO, they will result in the same distribution set and the same
decision will be chosen. However, these two distributions have very
different properties and the optimal decisions may be quite
different. Second, in the DRO approach, the worst-case distribution
for a decision is often unrealistic. For example, Scarf \cite{scarf}
shows that the worst-case distribution in the newsvendor context is
a two-point distribution. This raises the concern that the decision
chosen by this approach is guarding some overly conservative
scenarios, while performing poorly in more likely scenarios.
Unfortunately, these drawbacks seem to be inherent in the model
choice and cannot be satisfactorily remedied.

In this paper, we propose another choice of the distribution set in
the DRO framework that solves the above two drawbacks of the
mean-variance approach. Instead of using the mean and variance to
construct the distribution set, we choose to use the {\it likelihood
function}. More precisely, given a set of historical data, we define
the distribution set to be the set of distributions that make the
observed data achieve a certain level of likelihood. We call this
approach the {\it likelihood robust optimization} (LRO) approach.
The goal of this paper is to study the properties of LRO and its
performance.

First, we show that the LRO model is highly tractable. By applying
the duality theory, we formulate the robust counterpart of this
problem into a single convex optimization problem. In addition, we
show that our model is very flexible. We can add any convex
constraints (such as the moment constraints) to the distribution set
while still maintaining its tractability. Two concrete examples (a
newsvendor problem and a portfolio selection problem) are discussed
in the paper to illustrate the applicability of our framework.

Then we study the statistical theories behind the LRO approach by
illustrating the linkage between our approach and the Bayesian
statistics and empirical likelihood theory. We show that the
distribution set in our approach can be viewed as a confidence
region for the distributions given the set of observed data. Then we
discuss how to choose the parameter in the distribution set to
attain a specified confidence level. Furthermore, we show a
connection between the LRO approach and the mean-variance DRO
approach. Our analysis shows that the LRO approach is fully
data-driven, and it takes advantage of the full strength of the
available data while maintaining a certain level of robustness.

Finally, we test the performance of the LRO model in two problems, a
newsvendor problem and a portfolio selection problem. In the
newsvendor problem, we find that our approach produces similar
results compared to the mean-variance DRO approach when the
underlying distribution is symmetric, while the solution of our
approach is much better when the underlying distribution is
asymmetric. In the portfolio selection problem, we show by using
real historical data that our approach achieves decent returns.
Furthermore, the LRO approach will naturally diversify the
portfolio. As a result, the returns of the portfolio have a
relatively small fluctuation.

Following the initial version of this paper,\footnote{This paper is
based on the first author's tutorial paper in 2009.} Ben-Tal et al.
\cite{bental_divergence} studied a distributionally robust
optimization model where the distribution set is defined by
divergence measures. Their model contains the LRO as a special case.
They also discuss solvability and the statistical properties of
their models. In this paper, we focus on the distribution set
defined by the likelihood function and further explore the
connections to the empirical likelihood theory. In addition, we also
address the scenario when the sample space is continuous, which adds
insights to this class of approaches.

Besides the paper by Ben-Tal et al. \cite{bental_divergence}, a
recent paper by Bertsimas et al. \cite{bertsimas5} also studies DRO
with various choices of the distribution set. In particular, they
focus on using data and hypothesis-testing tools to construct those
sets. Although their paper is similar to ours, they do not use the
likelihood function to construct the robust distribution set nor do
they further study the properties of the DRO with such distribution
sets.

Two other papers related to this one are Iyengar \cite{Iyengar} and
Nilim and El Ghaoui \cite{Nilim}. In these two papers, the authors
study the robust Markov Decision Process (MDP) problem in which the
transition probabilities can be chosen from a certain set. They
mention the likelihood set as one choice. However, they do not
further explore the properties of this set or attempt to extend it
to general problems.

The remainder of this paper is organized as follows: In Section
\ref{sec:model}, we introduce our likelihood robust optimization
framework and discuss its tractability and statistical properties.
In Section \ref{sec:continuous}, we extend our discussions to
problems with continuous state space. In Section
\ref{sec:numerical}, we present numerical tests of our model.
Section \ref{sec:conclusion} concludes this paper.

\section{Likelihood Robust Optimization Model}
\label{sec:model}

In this section, we formulate the {\it likelihood robust
optimization} (LRO) model, discuss its tractability and statistical
properties. We only consider the case where the uncertainty
parameters have a finite discrete support in this section. We will
extend our discussions to the case with continuous support in
Section \ref{sec:continuous}.

Suppose we want to maximize an objective function $h(\bx,\bxi)$
where $\bx$ is the decision variable with feasible set $D$, and
$\bxi$ is a random variable taking values in $\Xi =
\{\bxi_1,\bxi_2,...,\bxi_n\}$. The set $\Xi$ is known in advance.
Assume we have observed $N$ independent samples of $\bxi$, with
$N_i$ occurrences of $\bxi_i$. We define:
\begin{eqnarray}\label{robustdistributionset}
{\mathbb D}(\gamma)=\left\{\bp=(p_1,p_2,...,p_n)\left|\sum_{i=1}^n
N_i\log{p}_i\ge\gamma, \sum_{i=1}^n p_i = 1, p_i\ge0, \forall
i=1,...,n\right.\right\}.
\end{eqnarray}
We call ${\mathbb D}(\gamma)$ the {\it likelihood robust
distribution set} with parameter $\gamma$. Note that ${\mathbb D}
(\gamma)$ contains all the distributions with support in $\Xi$ such
that the observed data achieves an empirical likelihood of at least
$\exp(\gamma)$. At this point, we treat $\gamma$ as a given
constant. Later we will discuss how to choose $\gamma$ such that
(\ref{robustdistributionset}) has a desirable statistical meaning.
We formulate the LRO problem as follows:
\begin{eqnarray}\label{lro}
\mbox{maximize}_{\bx\in D} && \left\{\min_{\bp\in{\mathbb
D(\gamma)}} \mbox{ } \sum_{i=1}^n p_ih(\bx,\bxi_i)\right\}.
\end{eqnarray}
In (\ref{lro}), we choose the decision variable $\bx$, such that the
expectation of the objective function under the worst-case
distribution is maximized, where the worst-case distribution is
chosen among the distributions such that the observed data achieve a
certain level of likelihood. Prior to this work, researchers have
chosen other types of distribution sets for the inner problem in
(\ref{lro}), e.g., distributions with moment constraints. Works of
that type have been reviewed in Section \ref{sec:intro}.

Here we comment on our assumption of the known support of the random
variable $\bxi$. In the LRO, the choice of $\Xi$ is important; a
different choice of $\Xi$ will result in a different distribution
set and a different solution to the decision problem. In practice,
sometimes $\Xi$ has a clear definition. For example, if $\bxi$ is
the transition indicator of a Markov chain (thus those $\bp$s are
transition probabilities), then the set of $\Xi$ is simply all the
attainable states in the Markov chain. However, in cases where the
choice of $\Xi$ is less clear, e.g., when $\bxi$ represents the
return of certain assets, the decision maker should choose $\Xi$ to
reflect his view of plausible outcomes. Also, as we will show later,
one can sometimes add other constraints such as moment constraints
into the distribution set. Once such constraints are added, the
choice of support often becomes less critical, since the support may
be constrained by existing constraints.

\subsection{Tractability of the LRO Model}
\label{subsec:tractability} In this subsection, we show that the LRO
model is easily solvable. To solve (\ref{lro}), we write down the
Lagrangian of the inner optimization problem:
\begin{eqnarray*}
L(\bp,\lambda,\mu) = \sum_{i=1}^n p_i
h(\bx,\bxi_i)+\lambda\left(\gamma-\sum_{i=1}^nN_i\log{p_i}\right)+\mu\left(1-\sum_{i=1}^np_i\right).
\end{eqnarray*}
Therefore, the dual formulation of the inner problem is
\begin{eqnarray}\label{dualinside}
\mbox{maximize}_{\lambda,\mu} &
\mu+\lambda\left(\gamma+N-\sum_{i=1}^nN_i\log{N_i}\right)-N\lambda\log{\lambda}+\lambda\sum_{i=1}^n
N_i\log{(h(\bx, \bxi_i) - \mu)} \nonumber\\
\mbox{s.t.} & \lambda \ge 0, \quad h(\bx,\bxi_i)-\mu\ge 0,
\quad\forall i.
\end{eqnarray}
Since the inner problem of (\ref{lro}) is always convex and has an
interior solution (given it is feasible), the optimization problem
(\ref{dualinside}) always has the same optimal value as (\ref{lro})
(see, e.g., Boyd and Vandenberghe \cite{boyd}). Then we combine
(\ref{dualinside}) and the outer maximization problem. We have that
(\ref{lro}) is equivalent to the following:
\begin{eqnarray}\label{wholeprogram}
\mbox{maximize}_{\bx,\lambda,\mu} & \mu
+\lambda(\gamma+N-\sum_{i=1}^nN_i\log{N_i})-N\lambda\log{\lambda}+\lambda\sum_{i=1}^nN_i\log{(h(\bx,\bxi_i)-\mu)}\nonumber\\
\mbox{s.t.} & \lambda \ge 0, \quad h(\bx,\bxi_i)-\mu\ge 0, \mbox{
}\forall i, \quad \bx\in D.
\end{eqnarray}
The following proposition follows immediately from examining the
convexity of the problem and studying its KKT conditions.
\begin{proposition}\label{solvability}
If $h(\bx,\bxi)$ is concave in $\bx$, then the likelihood robust
optimization problem described in (\ref{lro}) can be solved by
solving (\ref{wholeprogram}), which is a convex optimization
problem. Furthermore, if $(\bx^*,\lambda^*,\mu^*,\by^*)$ is an
optimal solution for (\ref{wholeprogram}), then
\begin{eqnarray*}
p_i=\frac{\lambda^* N_i}{h(\bx^*,\bxi_i)-\mu^*}
\end{eqnarray*}
is the corresponding worst-case distribution.
\end{proposition}


One advantage of the LRO approach is that one can integrate the
mean, variance and/or certain other information that is convex in
$\bp$ into the distribution set. For example, if one wants to impose
additional linear constraints on $\bp$, say $A\bp \ge \bb$ (note
that this includes moment constraints as a special case), then the
likelihood robust optimization model can be generalized as follows:
\begin{eqnarray}\label{generalizaion}
\mbox{maximize}_{\bx\in D} & \mbox{min}_{\bp\in {\mathbb D(\gamma)},
A\bp \ge \bb} & \sum_{i=1}^np_ih(x,\bxi_i).
\end{eqnarray}
By applying the duality theory again, we can transform
(\ref{generalizaion}) into the following problem ($\bmu$ is the dual
variable associated with the constraint $A\bp\ge \bb$, and we assume
that $\sum_{i=1}^n p_i = 1$ has been incorporated into $A\bp \ge
\bb$):
\begin{eqnarray*}
\mbox{maximize}_{\bx, \lambda, \bmu} & \bb^T\bmu
+(\gamma+N-\sum_{i=1}^nN_i\log{N_i})\cdot\lambda-N\lambda\log{\lambda}+\lambda\sum_{i=1}^nN_i\log{(h(\bx,\bxi_i)-\ba_i^T\bmu)}\\
\mbox{s.t.} & h(\bx,\bxi_i)-\ba_i^T\bmu\ge 0, \mbox{    }\forall i\\
& \bx\in D,\bmu\ge 0, \lambda \ge 0,
\end{eqnarray*}
which is again a convex program and readily solvable.

\subsection{Statistical Properties of the LRO Model}

In this subsection, we study the statistical theory behind the LRO
model. Specifically, we focus on the likelihood robust distribution
set defined in (\ref{robustdistributionset}). We will address the
following questions:
\begin{enumerate}
\item What are the statistical meanings of (\ref{robustdistributionset})?
\item How does one select a meaningful $\gamma$?
\item How does the LRO model relate to other types of distributionally
robust optimization models?
\end{enumerate}
We answer the first question in Section \ref{subsec:bayes} by using
the Bayesian statistics and empirical likelihood theory to interpret
the likelihood constraints. Those interpretations clarify the
statistical motivations of the LRO model. Then we answer the second
question in Section \ref{subsec:asymptotic} in which we perform an
asymptotic analysis of the likelihood region and point out an
asymptotic optimal choice of $\gamma$. We study the last question in
Section \ref{subsec:relationtoother} where we present a relationship
between our model and the traditional mean robust optimization
model.

\subsubsection{Bayesian Statistics Interpretation}
\label{subsec:bayes}

Consider a random variable taking values in a finite set $S$.
Without loss of generality, we assume $S =\{1,...,n\}$. Assume the
underlying probability distribution of the random variable is
$\bp=\{p_1,...,p_n\}$. We observe historical data
$\Psi=\{N_1,...,N_n\}$ in which $N_k$ represents the number of times
the random variable takes value $k$. Then the maximum likelihood
estimate (MLE) of $\bp$ is given by
$(N_1/N,N_2/N,...,N_n/N)$ where $N=\sum_{i=1}^n N_i$ is the total number of observations.

Now we examine the set of distributions $\bp$ such that the
likelihood of $\Psi$ under $\bp$ exceeds a certain threshold. We use
the concepts from Bayesian statistics (see e.g., Gelman et al.
\cite{gelman}). Instead of thinking that the data are randomly drawn
from the underlying distribution, we treat them as given, and we
define a random vector $\bp = (p_1,...,p_n)$ taking values on the
$(n-1)$-dimensional simplex $\Delta=\{\bp |p_1+\cdots
+p_n=1,p_i\geq0\}$ with probability density function proportional to
the likelihood function:
\begin{eqnarray*}
\prod_{i=1}^{n}{p_i}^{N_i}.
\end{eqnarray*}
This distribution of $\bp$ is known as the Dirichlet distribution. A
Dirichlet distribution with parameter $\Psi$ (denoted by
$Dir(\Psi)$) has the density function as follows:
\begin{eqnarray*}
f(\bp;\Psi)=\frac{1}{B(\Psi)}\prod_{i=1}^{n}{p_i}^{N_i-1}
\end{eqnarray*}
with
\begin{eqnarray*}
B(\Psi)=\int_{\bp\in\Delta} \prod_{i=1}^n{p_i}^{N_i-1}dp_1\cdots
dp_n=\frac{\prod_{i=1}^n\Gamma(N_i)}{\Gamma(\sum_{i=1}^n N_i)},
\end{eqnarray*}
where $\Gamma(\cdot)$ is the Gamma function. The Dirichlet
distribution is used to estimate the unknown parameters of a
discrete probability distribution given a collection of samples.
Intuitively, if the prior distribution is represented as
$Dir(\balpha)$, then $Dir(\balpha+\bbeta)$ is the posterior
distribution following a sequence of observations with histogram
$\bbeta$. For a detailed discussion on the Dirichlet distribution,
we refer the readers to \cite{gelman}.

In the LRO model, we assume a uniform prior $Dir(\be)$ on each point
in the support of the data, where $\be$ is the unit vector. After
observing the historical data $\Psi$, the posterior distribution
follows $Dir(\Psi+\be)$. Note that this process can be adaptive as
we observe new data.

Now we turn to examine the likelihood robust distribution set
\begin{eqnarray}\label{lroregion}
{\mathbb D}(\gamma) = \left\{\bp\left|\sum_{i=1}^n N_i\log{p_i} \ge
\gamma, \sum_{i=1}^n p_i = 1, p_i\ge 0, \forall i =
1,...,n\right.\right\}.
\end{eqnarray}
This set represents a region in the probability distribution space.
In the Bayesian statistics framework, we can compute the probability
that $\bp$ satisfies this constraint:
\begin{eqnarray*}\label{integral}
{\mathbb P}^*\left(\bp \in {\mathbb D}(\gamma)\right) =
\frac{1}{B(\Psi+\be)} \int_{\bp\in \Delta}
\prod_{i=1}^{n}{p_i}^{N_i}\cdot {\mathbb I}\left(\sum_{i=1}^n
N_i\log{p_i}\ge\gamma\right) dp_1\cdots dp_n,
\end{eqnarray*}
where ${\mathbb P}^*$ is the probability measure of $\bp$ given the
observed data $\Psi$  (thus ${\mathbb P}^*$ is the probability
measure under $Dir(\Psi + \be)$), and ${\mathbb I}(\cdot)$ is the
indicator function. When choosing the likelihood robust distribution
set, we want to choose $\gamma^*$ such that
\begin{eqnarray}\label{confidence}
{\mathbb P}^*\left(\bp \in {\mathbb D}(\gamma^*)\right) \ge 1-
\alpha
\end{eqnarray}
for some predetermined $\alpha$. That is, we want to choose
$\gamma^*$ such that (\ref{lroregion}) is the $1-\alpha$ confidence
region of the probability parameters. The LRO model can then be
interpreted as choosing the decision variable to maximize the
worst-case objective where the worst-case distribution is chosen
from the confidence region defined by the observed data.

However, in general, trying to find the exact $\gamma^*$ that
satisfies (\ref{confidence}) is computationally challenging. In the
next subsection, we study the asymptotic properties of $\gamma^*$,
which help us to approximate it.

\subsubsection{Asymptotic Analysis of the Likelihood Robust
Distribution Set} \label{subsec:asymptotic}

Now we investigate the asymptotic behavior of the likelihood robust
distribution set and give an explicit way to choose an appropriate
$\gamma$ in the LRO model. In this section, we assume that the true
underlying distribution of the data is $\bar{\bp} = \{\bar{p}_1,
...,\bar{p}_n\}$ with $\bar{p}_i
> 0$. We observe $N$ data drawn randomly from the underlying
distribution with $N_i$ observations on outcome $i$ ($\sum_{i=1}^n
N_i = N$). We define $\gamma_N$ to be the solution such that
\begin{eqnarray*}
\mathbb{P}_{\bN}\left(\bp \in {\mathbb D}(\gamma_N)\right) = 1 -
\alpha.
\end{eqnarray*}
Here $\alpha$ is a fixed constant and $\mathbb{P}_{\bN}$ is the
Dirichlet probability measure on $\bp$ with parameters
$N_1,...,N_n$. Clearly $\gamma_N$ depends on $N_i$ and thus is a
random variable. We have the following theorem about the asymptotic
properties of $\gamma_N$, whose proof is referred to Pardo
\cite{pardo}. (We use $\rightarrow_p$ to mean ``converge in
probability.'')

\begin{theorem}\label{thm:asymptotic_gamma1}
\begin{eqnarray*}\label{asymptotic1}
\gamma_N - \sum_{i=1}^n N_i\log{\frac{N_i}{N}} \rightarrow_p
-\frac{1}{2}\chi^2_{n-1, 1-\alpha}
\end{eqnarray*}
where $\chi^2_{d,1-\alpha}$ is the $1-\alpha$ quantile of a $\chi^2$
distribution with $d$ degrees of freedom.
\end{theorem}

Theorem \ref{thm:asymptotic_gamma1} provides a heuristic guideline
on how to choose the threshold $\gamma$ in the LRO model. In
particular, one should choose approximately
\begin{eqnarray*}\label{choiceofgamma}
\gamma^* = \sum_{i=1}^n {N_i}\log{\frac{N_i}{N}} -
\frac{1}{2}\chi^2_{n-1,1-\alpha}
\end{eqnarray*}
in order for the likelihood robust distribution set to have a
confidence level of $1-\alpha$. Note that the difference between
$\gamma^*$ and $\sum_{i} N_i\log{(N_i/N)}$ converges as $N$ grows.
This means that when the data set is large, the allowable
distributions must be very close to the empirical data. This is
unlike some other distributionally robust optimization approaches
which construct the distribution sets based on the mean and/or the
variance. Even with large data size, the distribution set in those
approaches may still contain distributions that are far from the
empirical distribution, such as a two-point distribution at the ends
of an interval. Therefore our approach does a better job in
utilizing the historical data to construct the robust region.

We also have a theorem about the convergence speed of $\gamma_N$. We
use $\Rightarrow_d$ to denote ``converge in distribution.''

\begin{theorem}\label{thm:asymptotic_gamma2}
\begin{eqnarray*}
\sqrt{N} \left(\frac{\gamma_N}{N} -
\sum_{i=1}^{n}\bar{p}_i\log{\bar{p}_i}\right) \Rightarrow_d X_0,
\end{eqnarray*}
where $X_0 = \sum_{i=1}^n(1+\log{\bar{p}_i})X_i$ with
$(X_1,...,X_n)\sim N(0, \Sigma)$ where
\begin{eqnarray*}
\Sigma = \left(\begin{array}{cccc} \bar{p}_1(1-\bar{p}_1) &
-\bar{p}_1\bar{p}_2 &
\cdots & -\bar{p}_1\bar{p}_n \\ -\bar{p}_1\bar{p}_2 & \bar{p}_2(1-\bar{p}_2) & \cdots & -\bar{p}_2\bar{p}_n \\
\vdots & \vdots & \ddots & \vdots \\ -\bar{p}_1\bar{p}_n &
-\bar{p}_2\bar{p}_n & \cdots &
\bar{p}_n(1-\bar{p}_n)\end{array}\right).
\end{eqnarray*}
\end{theorem}

{\noindent\bf Proof.} Given Theorem \ref{thm:asymptotic_gamma1}, it
suffices to prove that:
\begin{eqnarray*}\label{firstpart}
\sqrt{N}\left(\sum_{i=1}^n \frac{N_i}{N} \log{\frac{N_i}{N}} -
\sum_{i=1}^n \bar{p}_i\log{\bar{p}_i}\right) \Rightarrow_d X_0.
\end{eqnarray*}
To show this, we note that for fixed $N$, $(N_1,...,N_n)$ follows a
multinomial distribution with parameters $\bar{p}_1,...,\bar{p}_n$.
By Theorem 14.6 in Wassaman \cite{wassaman}, we have
\begin{eqnarray*}
\sqrt{N}\left(\frac{N_1}{N}-\bar{p}_1,...,\frac{N_n}{N} -
\bar{p}_n\right) \Rightarrow_d N(0,\Sigma).
\end{eqnarray*}
Define $f(x) = x\log{x}$. By the mean value theorem,
\begin{eqnarray*}
\sqrt{N}\left(\sum_{i=1}^n \frac{N_i}{N} \log{\frac{N_i}{N}} -
\sum_{i=1}^n \bar{p}_i\log{\bar{p}_i}\right) = \sum_{i=1}^n
f'(\eta_i)\sqrt{N}\left(\frac{N_i}{N}-\bar{p}_i\right)
\end{eqnarray*}
where $\eta_i$ lies between $N_i/N$ and $\bar{p}_i$. By the strong
law of large numbers, $\eta_i\rightarrow \bar{p}_i$ almost surely.
Therefore,
\begin{eqnarray*}
\sum_{i=1}^n f'(\eta_i)\sqrt{N}\left(\frac{N_i}{N}-\bar{p}_i\right)
\Rightarrow_d X_0,
\end{eqnarray*}
and thus the theorem is proved. \hfill $\Box$

\subsubsection{Relation to Other Types of DRO Models}
\label{subsec:relationtoother}

In this section, we show how one can relate other types of DRO
models to our LRO model through results in the empirical likelihood
theory. Given observed data $X_1,...,X_n$ with the empirical
distribution:
\begin{eqnarray*} \label{empiricaldistribution}
F_N(t) = \frac{1}{N} \sum_{i = 1}^N I(X_i\le t),
\end{eqnarray*}
we define the likelihood ratio of any distribution $F$ by
\begin{eqnarray*}\label{likelihoodratio}
R(F) = \frac{L(F)}{L(F_N)},
\end{eqnarray*}
where $L(F)$ is the likelihood value of the observations
$X_1,...,X_n$ under distribution $F$. Now suppose we are interested
in a certain parameter (or certain parameters) $\theta = T(F)$ of
the distribution. We can define the profile likelihood ratio
function as follows (see Owen \cite{owen}):
\begin{eqnarray}\label{profilerobust}
R_T(\theta) = \sup_F\{R(F)|T(F) = \theta\}.
\end{eqnarray}
Here some common choices of $T$ include the moments or quantiles of
a distribution. In DRO, one selects a region $\Theta$ and maximizes
the worst-case objective value for $\theta\in \Theta$. Using the
profile likelihood ratio function in (\ref{profilerobust}), we can
define $\Theta$ to be of the form
\begin{eqnarray*}
\Theta = \{\theta\left|R_T(\theta) \ge r_0\right.\}.
\end{eqnarray*}
When $\theta$ is the mean of the distribution, the next theorem
helps to determine $r_0$ such that asymptotically $\Theta$ is a
certain confidence interval for the true mean.
\begin{theorem}\label{thm:meanrobust}
{[Theorem 2.2 in Owen \cite{owen}]} Let $X_1,...,X_n$ be independent
random variables with common distribution $F_0$. Let $\mu_0={\mathbb
E}[X_1]$, and suppose that $0<Var(X_1)<\infty$. Then
$-2\log(R(\mu_0))$ converges in distribution to $\chi_{(1)}^2$ as
$n\rightarrow\infty$.
\end{theorem}

Therefore, in order for the set $\{\theta\left|R(\theta)\ge
r_0\right.\}$ to achieve a $1-\alpha$ confidence level (under the
Dirichlet distribution induced by the observed data), one should
approximately choose the boundary $\theta$ such that
\begin{eqnarray*}
R(\theta) = e^{-\frac{1}{2}\chi^2_{1,1-\alpha}}.
\end{eqnarray*}
Now we consider the following
optimization problem:
\begin{eqnarray}\label{find mu}
\mbox{maximize} & \sum_{i=1}^n N_i\log{p_i}\nonumber\\
\mbox{s.t.}& \sum_{i=1}^n p_i =1 \nonumber \\
&\sum_{i=1}^n d_ip_i=\theta\\
& p_i\geq 0, & i = 1,...,n.\nonumber
\end{eqnarray}
Let the optimal value of (\ref{find mu}) be $z(\theta)$. First,
$z(\theta)$ is concave in $\theta$ and achieves its maximum at
$\theta = \sum_{i=1}^n d_iN_i/N$ with the maximum value being
$\log(L(F_N))$. Also, when $\theta$ goes to $\max \{d_i|i =
1,...,n\}$ or $\min\{d_i|i = 1,...,n\}$, $z(\theta)$ goes to
$-\infty$. Therefore, there are exactly two $\theta$s such that
$\exp(z(\theta)) = L(F_N)\cdot\exp(-\chi^2_{1, 1-\alpha}/2)$. In
practice, we can find the two threshold $\theta$s by implementing a
bisection procedure. Denote the two threshold $\theta$s by $\theta_1
< \theta_2$. By Theorem \ref{thm:meanrobust}, we can construct a
mean robust optimization with uncertainty set $\theta\in [\theta_1,
\theta_2]$, and asymptotically this set has a probability of
$1-\alpha$ under the Dirichlet distribution defined by the observed
data.

Finally, we  note that the idea of using the profile likelihood
function can also be used to establish relationships between the LRO
model and other types of distributionally robust optimization.
However, there might not be a closed-form formula for the asymptotic
choice of $r_0$. In those cases, one may need to resort to sampling
methods to find an approximate one, the procedures of which are
beyond the discussion of this paper.

\section{Continuous State Space Case}
\label{sec:continuous}

In the previous section, we assume that the support of the uncertain
parameters is discrete. In the following, we extend our discussions
to the continuous case. We only consider the case where the
uncertain parameter is a scalar in this section.
%

It is tempting to directly extend the previous definition of the
likelihood region to the continuous case by using the probability
density function. However, with finite historical data, the
corresponding constraint will be defined only on a finite number of
points of the probability density function which effectively does
not make any restrictions to the distribution at all. Thus we have
to take a different approach. In this section, we propose an
approach that defines the robust distribution set by constructing a
band on the cumulative distribution function (CDF). Although
appearing to be different from our LRO model, this approach relies
on the same empirical likelihood theory discussed in the previous
section. We show that such an approach results in a tractable robust
counterpart and with proper choice of the band, the formulation is
statistically meaningful.

Given a set of observations $X_1,...,X_n$ drawn i.i.d. from an
underlying distribution $X$. A band on the CDF with support
$\{X_i\}$ and lower and upper bounds $\{L_i\}_{i=1}^n$ and
$\{U_i\}_{i=1}^n$ is defined by
\begin{eqnarray*}
\{F(\cdot)\mbox{ is a CDF} \left| L_i\le F(X_i)\le U_i,\mbox{ }
i=1,...,n\right.\}.
\end{eqnarray*}
Now we briefly discuss one example of such bands that is statistical
meaningful. We define the Kolmogorov-Smirnov band as:
\begin{eqnarray*}
\left\{F(\cdot)\left| \frac{i}{n} - D_{n, 1-\alpha}\le F(X_i) \le
\frac{i-1}{n} + D_{n,1-\alpha}, \mbox{ } i =1,...,n\right.\right\}.
\end{eqnarray*}
This band comes from the Kolmogorov-Smirnov goodness-of-fit test
with the statistics
\begin{eqnarray*}
D_n = \sup_{x}\left|F_n(x) - F(x)\right|,
\end{eqnarray*}
where $F_n(\cdot)$ is the empirical distribution of $X_1,...,X_n$.
By choosing $D_{n, 1-\alpha}$ such that $P(D_n \le
D_{n,1-\alpha})\ge 1-\alpha$, the band covers the true CDF
$1-\alpha$ fraction of times. This method can be modified to
construct a weighted Kolmogorov-Smirnov band where different weights
are used at different points. We refer the readers to Mason and
Schuenemeyer \cite{mason} for related discussions.

Once we have obtained such bands, we can write the corresponding
robust program as follows:
\begin{eqnarray}\label{dist2}
\mbox{maximize}_{\bx\in D} &\mbox{min}_F & {\mathbb E}_F[h(\bx,\xi)]\\
&\mbox{s.t.}& L_i\le F(X_i)\le U_i,\quad i = 1,...,n.\nonumber\\
&& F(\cdot) \mbox{ is a CDF. } \nonumber
\end{eqnarray}
By writing the CDF as the integral of the probability density
function (PDF), we can write (\ref{dist2}) as a semi-infinite
program:
\begin{eqnarray*}\label{density}
\mbox{maximize}_{\bx\in D} &\mbox{min}_f &\int_{-\infty}^{\infty}h(\bx,\xi)f(\xi)d\xi\nonumber\\
&\mbox{s.t.}&  L_i\le \int_{-\infty}^{X_i}f(\xi)d\xi\le U_i,\quad i=1,...,n\\
&& \int_{-\infty}^{\infty}f(\xi)d\xi =1, f(\xi)\ge 0.\nonumber
\end{eqnarray*}
By using the duality theorem, we can write the dual of the inner
program as follows:
\begin{equation}\label{temp3}
\begin{array}{rl}
\mbox{maximize}_{\by,\bz,\lambda}& \sum_{i=1}^n z_iL_i - \sum_{i=1}^n y_iU_i +\lambda\\
\mbox{s.t.}&
h(\bx,\xi)+\sum_{i=k}^{n}y_i-\sum_{i=k}^{n}z_i-\lambda\ge 0,\quad
 \forall \xi\in(X_{k-1},X_k],\quad k=1,...,n+1\\
& y_i\ge 0,z_i\ge 0, \forall i,
\end{array}
\end{equation}
where we define $X_0 = -\infty$ and $X_{n+1} = \infty$. If
$h(\bx,\xi)$ is concave in $\xi$, then the constraints in
(\ref{temp3}) can be reduced to:
\begin{eqnarray*}
h(\bx, X_{k-1})\ge\sum_{i=k}^{n}z_i - \sum_{i=k}^{n}y_i+\lambda, \quad\forall k\\
h(\bx, X_k)\ge\sum_{i=k}^{n}z_i-\sum_{i=k}^{n}y_i+\lambda,
\quad\forall k.
\end{eqnarray*}
Then combined with the outer problem, we obtain the robust
counterpart of (\ref{dist2}):
\begin{equation}\label{finalconinuoutsspace}
\begin{array}{rl}
\mbox{maximize}_{\bx,\by,\bz,\lambda} &\sum_{i=1}^n z_iL_i - \sum_{i=1}^n y_iU_i +\lambda\\
\mbox{s.t.}& h(\bx, X_{k-1})\ge\sum_{i=k}^{n}z_i - \sum_{i=k}^{n}y_i+\lambda, \quad\forall k\\
& h(\bx, X_k)\ge\sum_{i=k}^{n}z_i-\sum_{i=k}^{n}y_i+\lambda,
\quad\forall k \\
&y_i\ge 0,z_i\ge 0, \forall i, \mbox{  } \bx\in D.
\end{array}
\end{equation}
Further, if $h(\bx,\xi)$ is concave in $\bx$, then
(\ref{finalconinuoutsspace}) will be a convex program with a finite
number of variables and constraints, and thus can be solved easily.

\section{Applications and Numerical Results}
\label{sec:numerical}

In this section, we show two applications of the proposed LRO model
and perform numerical tests on them. The two applications are the
newsvendor problem (see Section \ref{subsec:newsvendor}) and the
portfolio selection problem (see Section \ref{subsec:portfolio}).

\subsection{Application 1: Newsvendor Problem}
\label{subsec:newsvendor}

In this subsection, we apply the LRO model to the newsvendor
problem. In such problems, a newsvendor facing an uncertain demand
has to decide how many newspapers to stock on a newsstand. If he
stocks too many, there will be a per-unit overage cost for each copy
that is left unsold; and if he stocks too few, there will be a
per-unit underage cost for each unmet demand. The problem is to
decide the optimal stocking quantity in order to minimize the
expected cost.

The newsvendor problem is one of the most fundamental problems in
inventory management and has been studied for more than a century.
We refer the readers to Khouja \cite{Khouja} for a comprehensive
review of this problem. In the classical newsvendor problem, one
assumes that the distribution of the demand is known (with
distribution $F$). The problem can then be formulated as:
\begin{eqnarray}\label{newsvendormodel}
\mbox{minimize}_{x} & G_F(x)=b{\mathbb E}_F(d-x)^++h{\mathbb
E}_F(x-d)^+.
\end{eqnarray}
In (\ref{newsvendormodel}), $x$ is the stocking quantity, $d$ is the
random demand with a probability distribution $F$, and $b$, $h > 0$
are the per-unit underage and overage costs, respectively.

It is well-known that a closed-form solution is available for this
problem. However, such a solution relies on the accurate information
of the demand distribution $F$. In practice, one does not always
have such information. In many cases, what one has is only some
historical data. To deal with this uncertainty, Scarf \cite{scarf}
proposed a distributionally robust approach in which a decision is
selected to minimize the expected cost under the worst-case
distribution, where the worst-case distribution is chosen from all
distributions that have the same mean and variance as the observed
data. However, as we discussed earlier, such an approach does not
fully utilize the data. It is also overly conservative as the
corresponding worst-case distribution for any given decision is a
two-point one, which is implausible from a practical point of view.


On the other hand, researchers have also proposed pure data-driven
models to solve the problem, i.e., using the empirical distribution
as the true distribution. However, pure data-driven approaches tend
to be less robust, since they ignore the potential deviations from
the data and do not guard against it.


In the following, we propose the LRO model for the newvendor
problem. We assume that the support of all possible demands is
$S=\{1,...,n\}$, i.e., the demand takes integer values between $1$
and $n$. Suppose we are given some historical demand data. We use
$N_i$ to denote the number of times that demand is equal to $i$ in
the observed historical data. The total number of observations is
$N=\sum_{i=1}^n N_i$.

The LRO model for the newsvendor problem is given as follows:
\begin{eqnarray}\label{lrofornewsvendor}
\mbox{minimize}_{x} & \mbox{max}_{\bp} & \sum_{i=1}^n p_i\left(b(d_i-x)^{+}+h(x-d_i)^{+}\right)\nonumber\\
&\mbox{s.t.} & \sum_{i=1}^n N_i \log{p_i} \ge \gamma \\
&& \sum_{i=1}^n p_i=1, \quad p_i\geq 0, \mbox{ }\forall i.\nonumber
\end{eqnarray}
Here the outer problem chooses a stocking quantity, while the inner
problem finds the worst-case distribution and cost with respect to
that decision. Here the first constraint makes sure that the
worst-case distribution is chosen among those distributions that
make the observed data achieve a chosen level of likelihood
$\gamma$. By applying the same techniques as in Section
\ref{sec:model}, we can write (\ref{lrofornewsvendor}) as a single
convex optimization problem:
\begin{eqnarray*}\label{finalnewsvendor}
\mbox{minimize}_{x, \lambda,\mu, \by} &\mu+ \lambda(\sum_{i=1}^n N_i\log{N_i}-N-\gamma)+ N\lambda\log{\lambda} - \sum_{i=1}^n N_i\lambda\log(y_i)\\
\mbox{s.t.} & b(d_i-x)+y_i\leq\mu \mbox{,     }\forall i\\
& h(x-d_i)+y_i\leq\mu \mbox{,     }\forall i\\
& \lambda\geq 0,\quad\by\geq 0.
\end{eqnarray*}
This problem is thus easily solvable. Similar approach can be
applied to multi-item newsvendor problems.

Now we perform numerical tests for this problem using the LRO model
and compare the solution to that of other methods. In the following,
we consider a newsvendor problem with unit underage and overage
costs $b=h=\$1$. We consider two underlying demand distributions.
The first one is a normal distribution with mean $\mu = 50$ and
standard deviation $\sigma = 50$. The second one is an exponential
distribution with $\lambda = 1/50$. We assume that both demand
distributions are truncated at 0 and 200.

For each underlying distribution, we perform the following
procedures for the LRO approach:
\begin{enumerate}
\item Generate $1000$ historical data from the underlying distribution, and record the number of times
the demand is equal to $i$ by $N_i$.
\item Construct the likelihood robust distribution set
${\mathbb D}(\gamma) = \{\bp|\sum_{i=1}^n N_i \log{p}_i \ge \gamma,
\sum_{i=1}^n p_i =1, p_i \ge 0, \forall i\}$. To choose a proper
$\gamma$, we use our asymptotic result in Theorem
\ref{thm:asymptotic_gamma1}. We choose $\gamma$ such that ${\mathbb
D}(\gamma)$ covers the true distribution with probability 0.95,
where the probability is under the Dirichlet measure with parameters
$N_1 + 1,...,N_n+ 1$. By Theorem \ref{thm:asymptotic_gamma1},
$\gamma$ should be approximately
$\sum_{i=1}^nN_i\log{\frac{N_i}{N}}-\chi^2_{100,0.95}$.
\item Solve the LRO model (\ref{lrofornewsvendor}) with the $\gamma$ chosen in Step
2.
\end{enumerate}
Using the same group of sample data, we also test the following
approaches:
\begin{enumerate}
\item LRO with fixed mean and variance:  Denote the sample mean by $\hat{\mu}$ and the sample variance by $\hat{\sigma}^2$. We add two constraints
$\sum_{i=1}^np_id_i = \hat{\mu}$ and $\sum_{i=1}^np_id_i^2 =
\hat{\mu}^2 + \hat{\sigma}^2$ to the inner problem of
(\ref{lrofornewsvendor}), that is, we only allow those distributions
that have the same mean and variance as the sample mean and
variance. As we discussed in (\ref{generalizaion}), this would still
be a convex problem and one can obtain the optimal solution easily.
We denote this approach by LRO($\hat{\mu},\hat{\sigma}^2$).
\item The distributionally robust approach with fixed mean and variance: This is the
approach proposed by Scarf \cite{scarf}. In this approach, one
minimizes the worst-case expected cost where the worst-case
distribution is chosen from all the distributions with a fixed mean
(equal to the sample mean) and a fixed variance (equal to the sample
variance). We call this approach the Scarf approach. In
\cite{scarf}, it is shown that the optimal decision for the Scarf
approach can be expressed in a closed-form formula:
\begin{eqnarray*}\label{scarfsolution}
x^*_{scarf} = \hat{\mu} +
\frac{\hat{\sigma}}{2}\left(\sqrt{\frac{b}{h}}-\sqrt{\frac{h}{b}}\right).
\end{eqnarray*}
\item The empirical distribution approach: We solve the optimal
solution using the empirical distribution as the true distribution.
\item We solve the optimal solution using the true underlying
distribution.
\end{enumerate}

\begin{table}[h]
\begin{center}
\begin{tabular}{|c|c|c|c|c|}
  \hline
                & Solution ($x^*$)  &  ${\mathbb E}_{F_{true}}[h(x^*,\xi)]$ & $h_{\mbox{LRO}}(x^*)$ \\ \hline
   LRO          & 56                &    35.97 (+0.84\%)           & 49.32 (+0,00\%) \\ \hline
   LRO($\hat{\mu},\hat{\sigma}^2$) & 54 &35.80 (+0.37\%)           & 49.47 (+0.30\%) \\ \hline
   Scarf        & 53                &    35.74 (+0.20\%)           & 49.63 (+0.63\%) \\ \hline
   Empirical    & 50                &    35.79 (+0.34\%)           & 50.03 (+1.40\%) \\ \hline
   Underlying   & 51                &    35.67 (+0.00\%)           & 49.87 (+1.11\%) \\ \hline
\end{tabular}\caption{Results when the underlying distribution is normal}\label{table:comparenormal}
\end{center}
\end{table}

\begin{table}[h]
\begin{center}
\begin{tabular}{|c|c|c|c|c|}
  \hline
                        & Solution ($x^*$)  &  ${\mathbb E}_{F_{true}}[h(x^*,\xi)]$ & $h_{\mbox{LRO}} (x^*)$ \\ \hline
   LRO              & 46              &  35.72(+3.51\%)   &   53.71 (+0.00\%)   \\ \hline
   LRO($\hat{\mu},\hat{\sigma}^2$)       & 41              &  34.99(+1.39\%)   &   54.45 (+1.38\%)  \\ \hline
   Scarf            & 50              &  36.64(+6.17\%)   &   54.27 (+1.04\%) \\ \hline
   Empirical        & 37              &  34.55(+0.16\%)   &   55.75 (+3.80\%)  \\ \hline
   Underlying       & 35              &  34.51(+0.00\%)   &   56.20 (+4.64\%)   \\ \hline
\end{tabular}\caption{Results when the underlying distribution is exponential}\label{table:compareexponential}
\end{center}
\end{table}

Our test results are shown in Tables \ref{table:comparenormal} and
\ref{table:compareexponential}. All the computations are run on a PC
with 1.80GHz CPU and Windows 7 Operating System. We use MATLAB
version R2010b to develop the algorithm for solving the optimization
problems. In Tables \ref{table:comparenormal} and
\ref{table:compareexponential}, the second column is the optimal
decision computed by each model. The third column shows the expected
cost of each decision under the true underlying distribution, which
is truncated $N(50,2500)$ in Table \ref{table:comparenormal} and
truncated $\mbox{Exp}(1/50)$ in Table
\ref{table:compareexponential}. The last column shows the objective
value of each decision under the LRO model, i.e., the worst-case
expected objective value when the distribution could be chosen from
${\mathbb D} (\gamma)$. The numbers in the parentheses in both
tables are the relative differences between the current solution and
the optimal solution under the measure specified in the
corresponding column.

\begin{figure}[h]
\centering \subfigure[Scarf model]{
\includegraphics[width=2.1in]{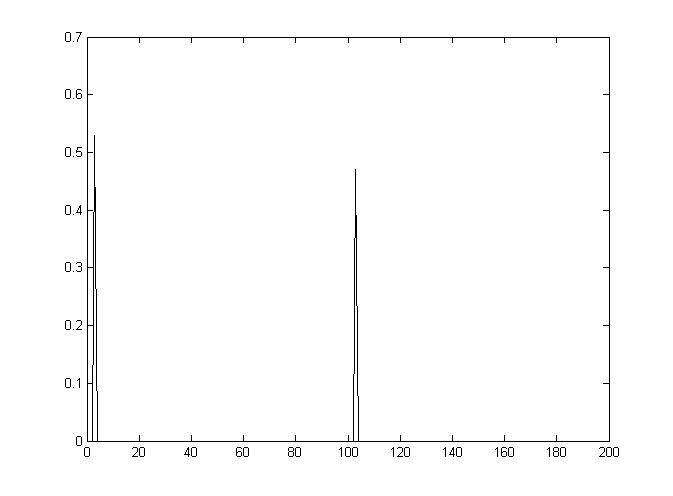}
\label{fig:scarf} }\subfigure[LRO model]{
\includegraphics[width=2.1in]{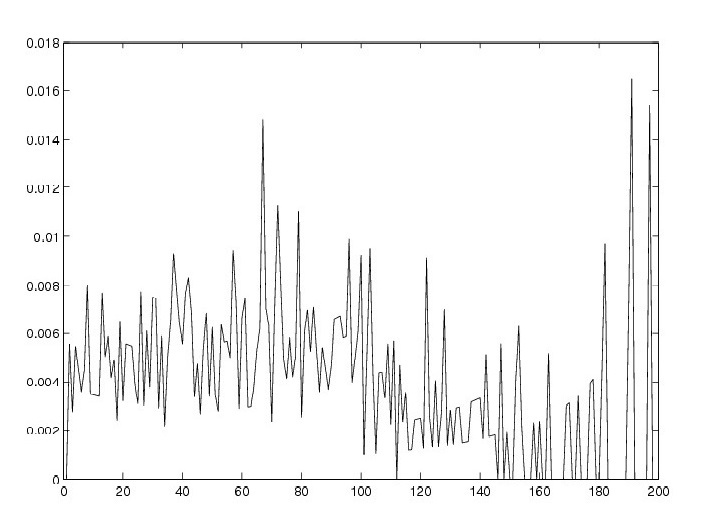}
\label{fig:lro} }\subfigure[LRO ($\hat{\mu},\hat{\sigma}^2$)
model]{\includegraphics[width = 2.1in]{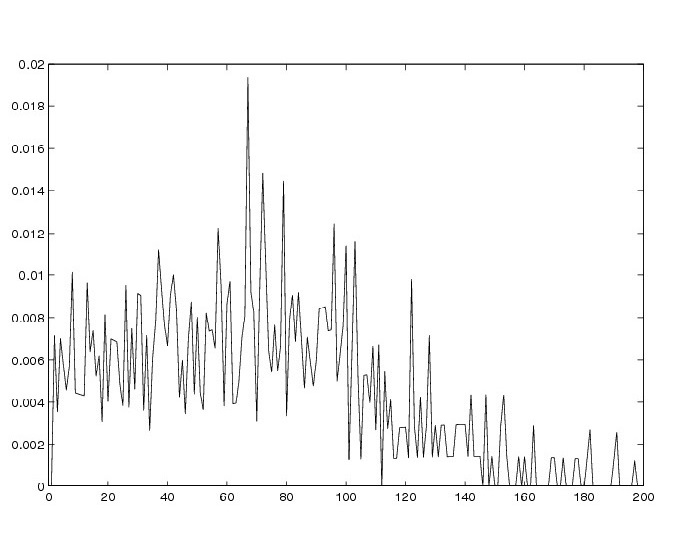}\label{fig:lrom} }
\caption{Worst-case distributions}
\label{fig:worstcasedistributions}
\end{figure}

In Table \ref{table:comparenormal}, one can see that when the
underlying demand distribution is normal, the solutions of different
models do not differ much from each other. This is mainly because
the symmetric property of the normal distribution pushes the
solution to the middle in either model. Nevertheless, the worst-case
distributions in different cases are quite varied. We plot them in
Figure \ref{fig:worstcasedistributions}.

In Figure \ref{fig:worstcasedistributions}, the $x$-axis is the
demand and the $y$-axis is the probability mass function over the
demands. We can see that the worst-case distribution in the Scarf
model is a two-point distribution (with positive mass at $3$ and
$103$). However, a two-point distribution is not a realistic one,
which means that Scarf's approach might be guarding some overly
conservative scenarios. In contrast, the worst-case distributions
when we use LRO and LRO($\hat{\mu},\hat{\sigma}^2$) are much closer
to the empirical data and look much more plausible. In particular,
the LRO with mean and variance constraints results in a worst-case
distribution closest to the empirical data among these three
approaches.

The situation is different in the second case when the underlying
distribution is an exponential distribution. In this case, Scarf's
solution is significantly worse than the LRO solutions. This is
because it does not use the information that the data are skewed. In
fact, it still only takes the mean and variance as the input. In
contrast, the LRO and LRO($\hat{\mu},\hat{\sigma}^2$) adapt to the
asymmetry of the data and only consider the distributions that are
close to the empirical distribution.

In both cases, we observe that the LRO($\hat{\mu},\hat{\sigma}^2$)
seems to work better. Indeed, we find that adding some constraints
on $\bp$ regarding the mean/variance of the distribution usually
helps the performance. This is because it helps to further
concentrate the distribution set to those that have similar shape as
the empirical distribution. Lastly, we find that the empirical
distribution works well if we use the true distribution to evaluate
the performance. However, it is not robust to potential deviations
to the underlying distribution. As shown in the last column in
Tables \ref{table:comparenormal} and \ref{table:compareexponential},
when we allow the underlying distribution to change to some degree
(within a 95\% confidence range), the performance of the solution
obtained by using the empirical distribution might be less than
optimal.

\subsection{Application 2: Portfolio Selection Problem}
\label{subsec:portfolio}

In this subsection, we apply the LRO model to the portfolio
selection problem. In the portfolio selection problem, there are $d$
assets that are available to invest for the next day. The decision
maker observes $N$ historical daily returns of the $d$ assets
$\{\bxi_i\}_{i=1}^N$ where each $\bxi_i$ is a $d$-dimensional vector
drawn from some underlying unknown distribution. In our case, we
consider a certain support $\Omega$ of all possible returns, and the
choice of $\Omega$ can be derived by using statistical models to
calculate the boundaries of the possible returns. The decision in
this problem is to make a portfolio selection $\bx$ in a feasible
set $D$ (e.g., $D$ may contain some budget constraints) to maximize
the expected return. We refer the readers to Luenberger
\cite{luenberger} for a thorough review of the literature on this
problem.

In this section, we demonstrate how to apply the LRO model to the
portfolio selection problem and compare its performance to that of
some other methods. First, we formulate the LRO model for the
portfolio selection problem as follows:
\begin{eqnarray}\label{lroforportfolio}
\mbox{maximize}_{\bx\in D} & \mbox{min}_{\bp} & \sum_{i=1}^N
p_i\cdot\bxi_i^T
\bx+\sum_{\sigma\in\Omega}p_\sigma\cdot\bxi_\sigma^T \bx\nonumber\\
& \mbox{s.t.} & \sum_{i=1}^N\log{p_i}\geq\gamma\\
&&\be^T\bp=1,\mbox{  }\bp \ge 0.\nonumber
\end{eqnarray}
In (\ref{lroforportfolio}), the first term in the objective function
corresponds to the scenarios that have been observed in the data,
while the second term corresponds to the unobserved scenarios (it
will be an infinite sum if $\Omega$ has an infinite support). The
first constraint is the likelihood constraint as we introduced in
Section \ref{sec:model}. For the ease of notation, we assume that
each return profile occurs only once in the data, which is
legitimate if we assume that the return is continuously distributed.
The choice of $\Omega$ can either be discrete or continuous. In
fact, as long as we can solve
\begin{eqnarray}\label{g}
g(\bx)=\min_{ \mbox{ }\bxi\in\Omega} \bxi^T\bx
\end{eqnarray}
efficiently for any $\bx$, the whole problem can be formulated as a
convex program with polynomial size. To see that, we take the dual
of the inner problem of (\ref{lroforportfolio}) and combine it with
the outer problem. We have that the entire problem can be written as
one single optimization problem as follows:
\begin{eqnarray}\label{finalportfolio}
\mbox{maximize}_{\bx,\mu,\lambda}&
\mu+\gamma\lambda+N\lambda-N\lambda\log{\lambda}+\lambda\sum_{i=1}^N\log(\bxi_i^T\bx-\mu)\nonumber\\
\mbox{s.t.}&\mu\leq\bxi_\sigma^T \bx,\quad\forall\sigma\in\Omega\\
&\lambda\geq 0, \quad \bx\in D.\nonumber
\end{eqnarray}
It appears that the first set of constraints could be of large size,
however, if $g(\bx)$ defined in (\ref{g}) can be solved efficiently,
then this is effectively one single constraint. For example, if
$\Omega=\{(r_1,...,r_d)|r_i\in[\underline{r}_i,\bar{r}_i]\}$, then
the constraint is equivalent to $\mu\leq\underline{\br}^T\bx$. If
$\Omega=\{(r_1,...,r_d)| ||\br-\br_0||_2\leq\eta\}$, then the
constraint is equivalent to $\mu\leq \br_0^T\bx-\eta||\bx||$, which
is also convex.

Again, we can add any convex constraints of $\bp$ into
(\ref{lroforportfolio}) such as constraints on the moments of the
return. Thus our model is quite flexible. Next we perform numerical
tests using our model. We adopt the following setup (this setup is
also used in other studies of the portfolio selection problem, see,
e.g., Delage and Ye \cite{delage}):
\begin{itemize}
\item We gather the historical data of 30 assets from the S\&P $500$ Index during the time period from 2001 to 2004. In each experiment, we choose four assets to focus on and
the decision is to construct a portfolio using these four assets for
each day during this period. We use the past 30 days' data as the
observed data in the LRO approach.
\end{itemize}

To select the $\gamma$ in (\ref{finalportfolio}), we still use the
asymptotic result in Theorem \ref{thm:asymptotic_gamma1}. In
particular, we choose the degree of freedom in the chi-square
distribution to equal the number of historical observations.
Although this is a heuristic choice, we find it works well in the
numerical test. We compare our approach to the following three other
approaches:
\begin{enumerate}
\item A naive approach in which each day, the stock with the highest
past 30 days' average daily return is chosen to be the sole stock in
the portfolio. We call this approach the {\it single stock} (SS)
approach.
\item An approach in which each stock is chosen with 1/4 weight (in terms of the capital
size) for each day. We call this approach the {\it equal weights}
(EQ) approach.
\item The distributionally robust approach with fixed mean and variance (see Popescu \cite{popescu}), called the DRO approach.
\end{enumerate}

We acknowledge that there are many other ways one can choose the
portfolio and we do not intend to show that our approach is the best
one for this problem. Instead, we aim to illustrate that the LRO
approach gives a decent performance with some desired features.

In Table \ref{table:portfolio}, we present the test results. The
results are obtained using codes written in MATLAB 2010b. Some of
the codes are borrowed from Delage and Ye \cite{delage}. We did
$100$ experiments (each with four different stocks randomly chosen
from the 30 to form the portfolio for each day), and the tested
period covers 721 days.
\begin{table}[h]
\centering
\begin{tabular}{||c|c|c|c||}
  \hline\hline
  Approach & \# LRO outperforms & Average gain of LRO over & Std of returns \\ \hline
  LRO   & N/A   & N/A     & $2.3\%$ \\
  SS    & 65    & $1.8\%$ & $3.3\%$ \\
  EQ    & 60    & $1.1\%$ & $1.5\%$ \\
  DRO   & 58    & $1.3\%$ & $2.4\%$ \\
  \hline\hline
\end{tabular}
\label{table:portfolio}\caption{Test results for the portfolio
selection problem}
\end{table}

In Table \ref{table:portfolio}, the numbers in the second column
indicate that among the $100$ experiments, the number of times the
overall return of the LRO approach outperforms the corresponding
method. The numbers in the third column are the improvements of the
average return of the LRO method over the corresponding model, and
the numbers in the last column are the standard deviations of all
daily returns of the portfolio constructed using each method. From
Table \ref{table:portfolio}, we observe that among the $100$
experiments, LRO outperforms all other tested methods in terms of
the number of times it has a higher return, as well as the average
overall return during the $721$ days. Also, it has a decent standard
deviation of the daily returns. In particular, it has a much smaller
standard deviation than the SS approach, and a comparable standard
deviation to the DRO approach. To investigate further, when we look
at the portfolios constructed by the LRO method, we find that $52\%$
of the time, it chooses a single stock to form the portfolio. And it
chooses two stocks, three stocks and all four stocks in its solution
for $39\%$, $8\%$ and $1\%$ of the times. Therefore, we find that
LRO approach implicitly achieves diversification, which is a feature
that is usually desired in such problems.

%
%


\section{Conclusion}
\label{sec:conclusion} In this paper, we propose a new
distributionally robust optimization framework, the likelihood
robust optimization model. The model optimizes the worst-case
expected value of a certain objective function where the worst-case
distribution is chosen from the set of distributions that make the
observed data achieve a certain level of likelihood. The proposed
model is easily solvable and has strong statistical meanings. It
avoids the over-conservatism of other robust models while protecting
the decisions from reasonable deviations from the empirical data. We
discuss two applications of our model. The numerical results show
that our model might be appealing in several applications.

\section{Acknowledgement}
The authors thank Dongdong Ge and Zhisu Zhu for valuable insights
and discussions. The authors also thank Erick Delage for insightful
discussions and for sharing useful codes for the numerical
experiments.

\bibliographystyle{plain}
\bibliography{lro}
\nocite{*}
\end{document}